\documentclass[12pt,a4paper]{amsart}

\usepackage{graphicx}
\usepackage[all]{xy}
\usepackage{comment}

\newcommand{\Z}{{\mathbb Z}}
\newcommand{\R}{{\mathbb R}}
\newcommand{\F}{{\mathbb F}}

\newtheorem{theorem}{Theorem}
\newtheorem{lemma}[theorem]{Lemma}

\numberwithin{equation}{section}

\begin{document}

\title{How to Solve a Diophantine Equation}

\author{Michael Stoll}

\address{Mathematisches Institut,
         Universit\"at Bayreuth,
         95440 Bayreuth, Germany.}
\email{Michael.Stoll@uni-bayreuth.de}

\date{23 February, 2010}

\begin{abstract}
  These notes represent an extended version of a talk I gave for
  the participants of the IMO~2009 and other interested people. We introduce
  diophantine equations and show evidence that it can be hard to solve them.
  Then we demonstrate how one can solve a specific equation related to numbers
  occurring several times in Pascal's Triangle with state-of-the-art methods.
\end{abstract}

\maketitle


\section{Diophantine Equations}
\label{sec:1}

The topic of this text is \emph{Diophantine Equations}. A diophantine equation
is an equation of the form
\[ F(x_1, x_2, \dots, x_n) = 0 \,, \]
where $F$ is a polynomial with integer coefficients, and one asks for
solutions in \emph{integers} (or rational numbers, depending on the problem).
They are named after Diophantos of Alexandria on whom not much is known
with any certainty. Most likely he lived around 300~AD. He wrote the
\emph{Arithmetika}, a text consisting of 13~books, a number of which have
been preserved. In this text, he explains through many examples ways of
solving certain kinds of equations like the above in rational numbers.
Diophantos was also one of the first to introduce symbolic notation for
the powers of an indeterminate.

To give you a flavor of this kind of question, let me show you some examples.
Ideally, you should cover up the part of the page below the equation and
try to find a solution for yourself before you read on. The first equation is
\[ x^3 + y^3 + z^3 = 29 \,, \]
an equation in three unknowns, to be solved in 
(not necessarily positive) integers.
I trust it did not take you very long to come up with a solution like
$(x,y,z) = (3,1,1)$ or maybe $(4,-3,-2)$. Now let us look at
\[ x^3 + y^3 + z^3 = 30 \,. \]
Try to solve it for a while before you look up a solution in this
footnote\footnote%
{$x = 2 \,220 \,422 \,932, \quad y = -2 \,218 \,888 \,517, \quad z = -283 \,059\, 965$.}.
This solution is the smallest and was found by computer search in July~1999
and published in~2007~\cite{Becketal}. This already indicates that it may
be quite hard to find a solution to a given diophantine equation.
Now consider
\[ x^3 + y^3 + z^3 = 31 \,. \]
Did you try to solve it? You should have come to the conclusion that there
is no solution: the third power of an integer is always $\equiv -1,0$ or
$1 \bmod 9$, so a sum of three cubes can never be $\equiv 4$ or $5 \bmod 9$.
Since $31 \equiv 4 \bmod 9$, the number $31$ cannot be a sum of three cubes. If we
replace $31$ with~$32$, the same argument applies. So we consider
\[ x^3 + y^3 + z^3 = 33 \]
next. If you were able to solve this, you should consider making diophantine
equations your research area. The sad state of affairs is that it is an
open problem whether this equation has a solution in integers or not!%
\footnote{This introduction was inspired by a talk Bjorn Poonen gave
at a workshop in Warwick in~2008.}

So the following looks like an interesting problem:
to decide if a given diophantine
equation is solvable or not. In fact, this problem appears on the most
famous list of mathematical problems, namely the 23~problems David Hilbert
stated in his address to the International Congress of Mathematicians
in Paris in~1900 as questions worth working on in the new century. The
description of the tenth problem in Hilbert's list reads thus
(in the German original~\cite{HilbertD}, see~\cite{HilbertE} for an English
translation of Hilbert's address):

\begin{center}
  \framebox{\includegraphics[width=0.98\textwidth]{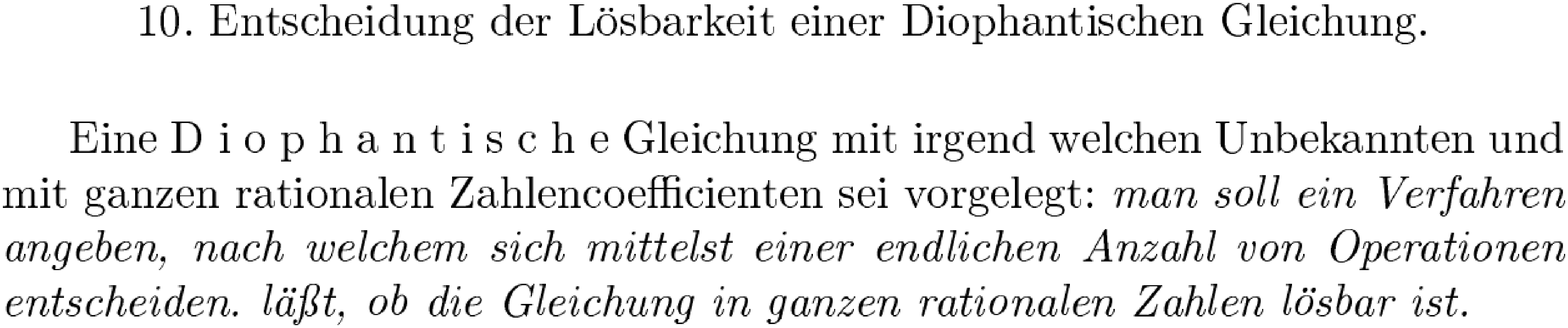}}
\end{center}

Here is an English translation.
\begin{quotation}
  Given a diophantine equation with any number of unknown quantities and
  with rational integral numerical coefficients: \emph{to devise a process according
  to which it can be determined by a finite number of operations whether the
  equation is solvable in rational integers.}
\end{quotation}

In modern terminology, Hilbert asks for an \emph{algorithm} that, given
a polynomial $F(x_1, \dots, x_n)$ with integral coefficients, decides whether
the equation
\[ F(x_1, \dots, x_n) = 0 \]
can be solved in integers. This is
commonly known as \emph{Hilbert's Tenth Problem}. It is not only the
shortest problem on Hilbert's list, it is also the only
decision problem\footnote{A \emph{decision problem} asks for an algorithm
that decides if a given element of a specified set has a specified property.},
so it is somewhat special. From the wording it can be inferred that Hilbert
believed in a positive solution to his problem: such an algorithm had to
exist. In fact, at the end of the introductory part of his speech, before
turning to the list of problems, he says 
\begin{quotation}
   \dots in der Mathematik giebt es kein Ignorabimus!
\end{quotation}
(There is no `Ignorabimus'\footnote{This Latin word means `we will not know'.}
in mathematics.) This indicates that Hilbert was convinced that
every mathematical problem must have a definite solution.

The simple examples I have shown at the beginning may (or should) have given
you a feeling that this problem may actually be very hard. This is also
what happened historically. People got more and more convinced that the
answer to Hilbert's Tenth Problem was likely to be negative: an algorithm
conforming to the given specification does not exist. Now if an algorithm 
does exist that performs a certain task, it is fairly clear how one can
prove this fact. Namely, one has to find such an algorithm and write it down,
then everybody will agree that it
indeed \emph{is} an algorithm solving the given problem. To show that such an
algorithm \emph{does not} exist is a quite different matter. One needs
some way of getting a handle on all possible algorithms, so that one can
show that none of them solves the problem. The relevant theory, which is
a branch of mathematical logic, did not yet exist when Hilbert gave his talk.
It was developed a few decades later, leading to such famous results as
G\"odel's Incompleteness Theorem, which definitely showed that there certainly
is an \emph{Ignorabimus} in mathematics. Indeed, work of several people,
most notably Martin Davis, Hilary Putnam and Julia Robinson, made it possible for
Yuri Matiyasevich to finally prove in~1970 the following result\footnote{%
See~\cite{Matiyasevich} for an accessible account of the problem and its solution}.

\begin{theorem}[Davis, Putnam, Robinson; Matiyasevich] \label{thm:1} \strut \\
   The solvability of diophantine equations is undecidable.
\end{theorem}

In fact, he proved a much stronger result, which implies for example
that there is an explicit polynomial $F(x_0,x_1,\dots,x_n)$ such that
there is no algorithm that, given $a \in \Z$ as input, decides whether
or not there is an integral solution to
\[ F(a, x_1, \dots, x_n) = 0 \,. \]
Note that if a diophantine equation is solvable, then we can prove it,
since we will eventually find a solution by searching through the
countably many possibilities (but we do not know beforehand how far we have
to search). So the really hard problem is to prove
that there are no solutions when this is the case. A similar problem
arises when there are finitely many solutions and we want to find them
all. In this situation one expects the solutions to be fairly small\footnote{%
The large solution to $x^3 + y^3 + z^3 = 30$ is no counterexample to this
statement, since there should be infinitely many solutions in this case.}.
So usually it is not so hard to find all solutions; what is difficult
is to show that there are no others.

So, given Theorem~\ref{thm:1}, should we give up all attempts to solve
diophantine equations, convinced that the task
is completely hopeless? That would be premature. We might still be able
to prove positive results when we restrict the set of equations in some
way. For example, there are quite good reasons to believe that
there should be a positive answer to Hilbert's question for equations
\emph{in two variables}. In the remainder of this contribution, we will
consider one such equation as an example case and show with what kind
of methods it can be attacked and solved.


\section{The Example Equation}
\label{sec:2}

\begin{figure}[htb]
   \includegraphics[width=\textwidth]{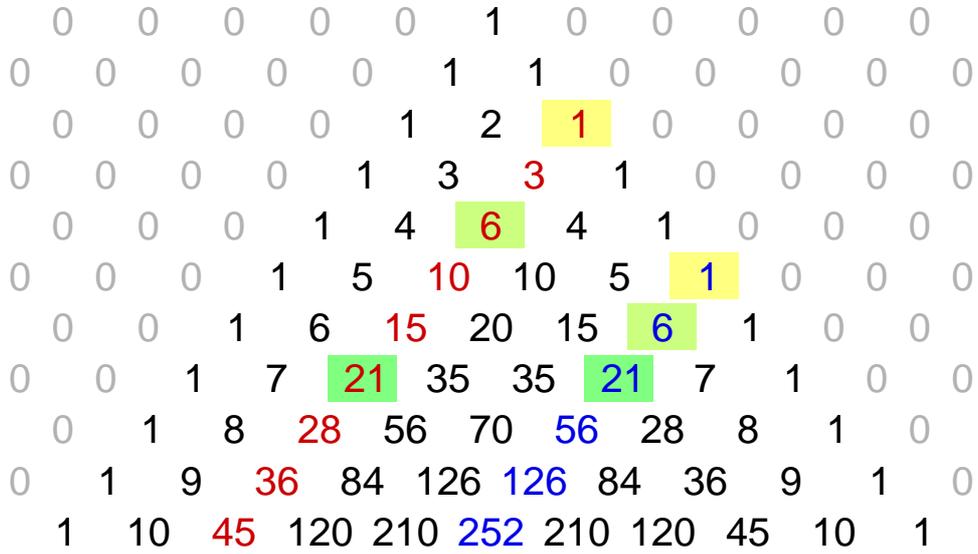}
   \caption{Pascal's Triangle}
   \label{fig:Pascal}
\end{figure}

The equation we want to consider here is motivated by the following
question. Consider Pascal's Triangle (Fig.~\ref{fig:Pascal}).
Which natural numbers occur several times in this triangle, if we
disregard the outer two ``layers'' ($1,1,1,\dots$ and $1,2,3,\dots$)
on either side and the obvious reflectional symmetry?

In other words, what are the integral solutions to the equation
\begin{equation} \label{bineq}
   \binom{y}{k} = \binom{x}{l} \,,
\end{equation}
subject to the conditions $1 < k \le y/2$, $1 < l \le x/2$ and $k < l$?
The following solutions are known.
\begin{eqnarray*}
   && \binom{16}{2}=\binom{10}{3},\quad 
\binom{56}{2}=\binom{22}{3},\quad 
\binom{120}{2}=\binom{36}{3},\\
   && \binom{21}{2}=\binom{10}{4},\quad 
\binom{153}{2}=\binom{19}{5},\quad 
\binom{78}{2}=\binom{15}{5}=\binom{14}{6},\\
   && \binom{221}{2}=\binom{17}{8},\quad 
\binom{F_{2i+2}F_{2i+3}}{F_{2i}F_{2i+3}}=\binom{F_{2i+2}F_{2i+3}-1}{F_{2i}F_{2i+3}+1} 
\mbox{ for } i=1,2,\ldots,
\end{eqnarray*}
where $F_n$ is the $n$-th Fibonacci number.

Equation~(\ref{bineq}) is not a diophantine 
equation according to our definition,
since it depends on $k$ and~$m$ in a non-polynomial way. Also, it is way
too hard to solve. So we specialize by fixing $k$ and~$l$. The cases
\[ (k,l) \in \{ (2, 3), (2, 4), (2, 6), (2, 8), (3, 4), (3, 6), (4, 6), (4, 8) \} \]
have already been solved completely, see~\cite{SdW}. Each of these cases requires some
deep mathematics of a flavor similar to what is described below.
The next interesting case is obviously $(k,l) = (2,5)$, leading to the equation
\begin{equation} \label{eqn:example}
   \binom{y}{2} = \binom{x}{5}\,, \quad\mbox{or}\quad
    60 y (y-1) = x (x-1) (x-2) (x-3) (x-4) \,.
\end{equation}
So we are asking for numbers that occur both in the red and the blue
diagonal in Figure~\ref{fig:Pascal}.

The first step in solving an equation like~(\ref{eqn:example}) is
to go and look for its solutions. We easily find solutions with
\[ x = 0, 1, 2, 3, 4, 5, 6, 7, 15 \quad\mbox{and}\quad 19 \,, \]
and then no further ones. (Only the last two are `nontrivial' in the
sense that they satisfy the constraints given above. Also, there are no
solutions with $x < 0$, since then the right hand side is negative,
but the left hand side can never be negative for $y \in \Z$.) This now raises
the question if we have already found them all, and if so, how to prove it.

This is a good point to look at what is known about the solution set of
equations like~(\ref{eqn:example}) in general. The first important result
was proved by Carl Ludwig Siegel in~1929. (See~\cite[Section~D.9]{HS} for
a proof.)

\begin{theorem}[Siegel] \label{thm:Siegel} \strut \\
   Let $F \in \Z[x,y]$. If the solutions to $F(x,y) = 0$ cannot be rationally
   parameterized, then $F(x, y) = 0$ has only 
finitely many solutions in integers.
\end{theorem}

A \emph{rational parameterization} of $F(x,y) = 0$ is a pair of rational
functions $f(t)$, $g(t)$ (quotients of polynomials), not both constant, such
that $F(f(t), g(t)) = 0$ (as a function of~$t$). The existence of such a
rational parameterization can be algorithmically checked; for our equation
it turns out that it is not rationally parameterizable. So we already know
that there are only finitely many solutions. In particular, we have a chance
that our list is complete. On the other hand, Theorem~\ref{thm:Siegel} and its
proof are inherently \emph{ineffective}: we do not get a bound on the size
of the solutions, so this result gives us no way of checking that our list
is complete. This somewhat unsatisfactory state of affairs did not change
until the 1960s, when Alan Baker developed his theory of `linear forms in
logarithms' that for the first time provided explicit bounds for solutions
of many types of equations. For this breakthrough, he received the Fields
Medal. Baker's results cover a class of equations that contains our
equation~(\ref{eqn:example}). For our case, what he proved comes down
to roughly the following.
\begin{equation} \label{eqn:baker}
   |x| < 10^{10^{10^{10^{600}}}} \,.
\end{equation}
This reduces the solution of our equation~(\ref{eqn:example}) to a finite
problem. The
inequality in~(\ref{eqn:baker}) gives us an explicit upper bound for~$x$.
So we only have to check
the finitely many possibilities that remain, and we will obtain the complete
set of solutions to~(\ref{eqn:example}). From a very pure mathematics viewpoint,
we may therefore consider our problem as solved.
On the other hand, from a more
practical point of view, we would like to actually obtain the complete
list of solutions, and the assertion that it is possible in principle to
get it does not satisfy us. To say that the number showing up
in~(\ref{eqn:baker}) defies all imagination is a horrible understatement,
and one cannot even begin to figure out how long it would take to actually
perform all the necessary computations.

However, time did not stop in the 1960s, and with basically still the same
method, but with a lot of refinements and improvements thrown in, we are
now able to prove the following estimate.
\begin{equation} \label{eqn:better}
   |x| < 10^{10^{600}} \,.
\end{equation}
You may rightly ask whether something has really been gained, in practical
terms. The number of electrons in the universe is 
estimated to be about~$10^{80}$,
so we cannot even write down a number with something like $10^{600}$~digits!
However, it will turn out that the improvement 
represented by~(\ref{eqn:better})
is crucial. But before we can see this, we need to look at our problem
from a different angle.


\section{A Geometric Interpretation}
\label{sec:3}

\begin{figure}[t]
   \includegraphics[width=0.625\textwidth]{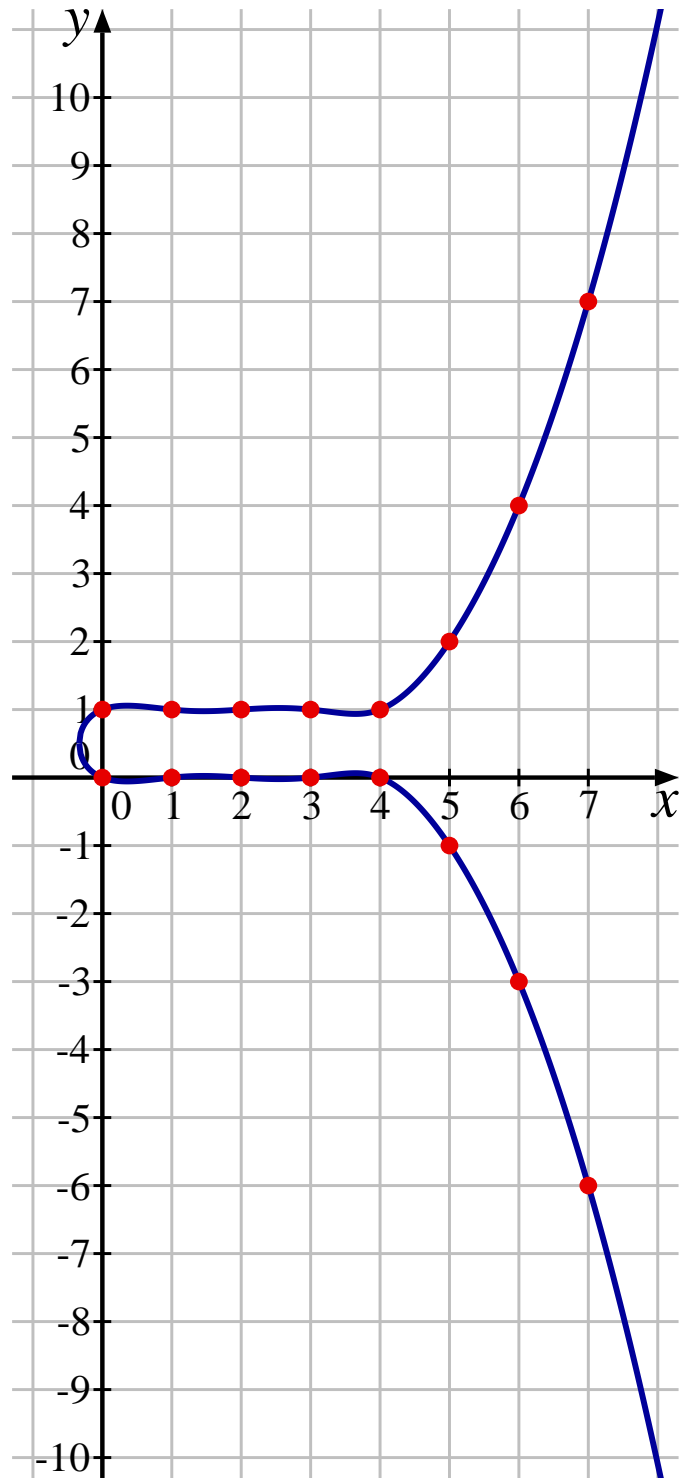}
   \caption{The curve given by~(\ref{eqn:example}), with some integral points.}
   \label{fig:curve}
\end{figure}

The idea is to translate our at first sight algebraic problem
((\ref{eqn:example}) is an algebraic equation) into a geometric one.
An equation $F(x,y) = 0$ in two variables defines a subset of the plane,
consisting of those points whose coordinates satisfy the equation. If
$F$ is a (non-constant) polynomial, this solution set is called a
\emph{plane algebraic curve}. We can draw the curve~$C$ corresponding
to our equation~(\ref{eqn:example}) in the real
plane~$\R^2$, see Fig.~\ref{fig:curve}. We are now interested
in the \emph{integral points} on~$C$, since they correspond to integral
solutions to~(\ref{eqn:example}). The set of integral points on~$C$ is
denoted~$C(\Z)$.

This set $C(\Z)$ of integral points on the curve~$C$ by itself does not
have any useful additional structure. But we can make use of a well-developed
theory, called Algebraic Geometry, that studies sets defined by a collection
of polynomial equations, and in particular algebraic curves like~$C$.
This theory tells us that we can \emph{embed} the curve~$C$ into another
object~$J$, which is not a curve, but a surface. This can be constructed
for any curve and is called the \emph{Jacobian variety} of the curve%
\footnote{The Jacobian variety need not be a surface; its dimension depends
on the curve.}.
The interesting fact about~$J$ (and Jacobian varieties in general) is that $J$
is a~\emph{group}. More precisely, there is a composition law on~$J$ that
is defined in a geometric way and that turns (for example) the set $J(\Z)$
of integral points\footnote{Algebraic geometers use the set of \emph{rational}
points here. This does not make a difference, since $J$ is a projective variety.
(Which means that the coordinates can be scaled so as to remove denominators.)}
on~$J$ into an \emph{abelian group.} In a similar spirit as Siegel's
Theorem~\ref{thm:Siegel} (and actually preceding it), we have the following
important result, valid for Jacobian varieties in general.
(See~\cite[Part~C]{HS}.)

\begin{theorem}[Weil 1928] \label{thm:Weil} \strut \\
   If $J$ is the Jacobian variety of a curve, then the abelian group $J(\Z)$
   is finitely generated.
\end{theorem}

This means that we can (in principle) get an explicit description of the
group~$J(\Z)$ in terms of generators and relations. If we have that, we
may be able to use the group structure and the geometry in some way to get
a handle on the elements of~$J(\Z)$ that are in the image of~$C$; these
correspond exactly to the integral points on~$C$.

In general, it is not known whether it is always possible to actually determine
explicit generators of a group like~$J(\Z)$ in an algorithmic way, although
there are some `standard conjectures' whose truth would imply a positive
answer. There are methods available that, with some luck, can find a set
of generators, but they are not guaranteed to work in all cases. This is
the point where the method we are describing may fail in practice. In our
specific example, we are lucky, and we can show that $J(\Z)$ is a
\emph{free abelian group of rank~6}:
\begin{equation} \label{eqn:gens}
   J(\Z) = \Z \, P_1 + \Z \, P_2 + \Z \, P_3 + \Z \, P_4 + \Z \, P_5 + \Z \, P_6
\end{equation}
with explicitly known points $P_1, \dots, P_6 \in J(\Z)$.

Let $\iota : C \to J$ denote the embedding of~$C$ into~$J$. The surface~$J$
lives in some high-dimensional space, and we can 
specify integral points on it by
a bunch of coordinates. We can measure the size of such a point by taking
the logarithm of the largest absolute value of the coordinates (this tells
us roughly how much space we need to write the point down). This gives us
a function
\[ h : J(\Z) \to \R_{\ge 0} \]
called the \emph{height}. One can show that this height function has the
following properties. The first one tells us how the height relates to the size of
integral points on our curve.
\begin{equation} \label{eqn:heightC}
   h\bigl(\iota(x,y)\bigr) \approx \log |x| 
\end{equation}
for points $(x,y) \in C(\Z)$ such that $x$ is not very small.

The second property says that the height function behaves well with
respect to the group structure on~$J$.
\begin{equation} \label{eqn:heightJ}
   h(n_1 P_1 + n_2 P_2 + n_3 P_3 + n_4 P_4 + n_5 P_5 + n_6 P_6)
     \approx n_1^2 + n_2^2 + n_3^2 + n_4^2 + n_5^2 + n_6^2 \,.
\end{equation}
(To be precise, each side can be bounded by an explicit constant multiple
of the other one. To be more precise, $h$ is, up to a bounded error, a
positive definite quadratic form on~$J(\Z)$.)

If we now combine the estimate~(\ref{eqn:better}) with the properties
(\ref{eqn:heightC}) and~(\ref{eqn:heightJ}) of the height~$h$, then we
obtain the following statement.

\begin{lemma}
   If $(x,y) \in C(\Z)$, then we have
   \[ \iota(x,y) = n_1 P_1 + n_2 P_2 + n_3 P_3 + n_4 P_4 + n_5 P_5 + n_6 P_6 \]
   with coefficients $n_j \in \Z$ satisfying $|n_j| < 10^{300}$.
\end{lemma}

Of course, the bound $10^{300}$ given here is not precise; a precise bound
can be given and is of the same order of magnitude.

The conclusion is that using the additional structure we have on~$J$ enables
us to reduce the size of the search space from about $10^{10^{600}}$ to `only'
$10^{1800}$ (there are six coefficients $n_j$ with about~$10^{300}$ possible
values each). This is, of course, still much too 
large to check each possibility
(think of the electrons in the universe), but, and this is the decisive point,
the numbers $n_j$ we have to deal with can be represented easily on a computer,
and we \emph{can compute} with them!


\section{Needles in a Haystack}
\label{sec:4}

We now have an enormous haystack
\[ H = \left\{(n_1, n_2, n_3, n_4, n_5, n_6) \in \Z^6 : |n_j| < 10^{300}\right\} \]
of about $10^{1800}$ pieces of grass that
contains a small number of needles. We want to find the needles. Instead of
looking at each blade of grass in the haystack, we can try to solve this
problem faster by finding conditions on the possible positions of the needles
that rule out large parts of the haystack. This is the point where we make
use of the fact that the group law on~$J$ is defined via geometry. Our objects
$C$, $J$ and~$\iota$ are defined over~$\Z$, 
therefore it makes sense to consider
their defining equations modulo~$p$ for prime numbers~$p$. We denote the
field $\Z/p\Z$ of $p$~elements by~$\F_p$. The sets of points with coordinates
in~$\F_p$ that satisfy these defining equations mod~$p$ are denoted by
$C(\F_p)$ and~$J(\F_p)$. Then for all but finitely many~$p$ (and the exceptions
can be found explicitly), $J(\F_p)$ is again an abelian group, and it
contains the image $\iota(C(\F_p))$ of~$C(\F_p)$. The group $J(\F_p)$ is
finite, and so is the set~$C(\F_p)$; both can be computed. Furthermore,
the following diagram commutes, and the geometric nature of the group structure
implies that the right hand vertical map is a group homomorphism.
\[ \xymatrix{ C(\Z) \ar[r]^{\iota} \ar[d] & J(\Z) \ar@{=}[r] \ar[d]
                    & \mbox{\raisebox{3pt}{$\Z^6$}} \ar[dl]^{\alpha_p} \\
               C(\F_p) \ar[r]^{\iota_p} & J(\F_p) &
             }
\]
The vertical maps are obtained by reducing the coordinates of the points
mod~$p$. The diagonal map~$\alpha_p$ is again a 
group homomorphism, determined by
the image of the generators $P_1, \dots, P_6$ of~$J(\Z)$. The following
is now clear.

\begin{lemma}
   Let $(x,y) \in C(\Z)$ and $\iota(x,y) = n_1 P_1 + \cdots + n_6 P_6$.
   Then
   \[ \alpha_p(n_1, n_2, n_3, n_4, n_5, n_6) \in 
\iota_p\bigl(C(\F_p)\bigr) \,. \]
\end{lemma}

The subset
$\Lambda_p = \alpha_p^{-1}\bigl(\iota_p(C(\F_p))\bigr) \subset \Z^6$ is
(usually, when $\alpha_p$ is surjective) a union of $\#C(\F_p)$ cosets
of a subgroup of index~$\#J(\F_p)$ in~$\Z^6$. Since one can show that
$\#C(\F_p) \approx p$
and $\#J(\F_p) \approx p^2$ (reflecting dimensions $1$ and~$2$, respectively),
we see that the intersection of our haystack~$H$ with~$\Lambda_p$ has only
about $1/p$~times as many elements as the original haystack.
This does not yet help very much, but we can try to \emph{combine} information
coming from \emph{many} primes. If $S$ is a (finite, but large) set of
prime numbers, then we set
\[ \Lambda_S = \bigcap_{p \in S} \Lambda_p
    \quad\mbox{and obtain}\quad
    \iota\bigl(C(\Z)\bigr) \subset \Lambda_S \cap H \,.
\]
If we make $S$ sufficiently large (about a thousand primes, say), then it
is likely that the set on the right hand side is quite small, so that we can
easily check the remaining possibilities. The idea is that the reductions
of the haystack size we obtain from several distinct primes should
accumulate, so that we can expect a reduction by a factor
which is roughly the product of all the primes in~$S$.

We have to be careful to select
the primes in a good way so that the description of the sets $\Lambda_S$ we
encounter on the way stays within a reasonable complexity. It is, however,
indeed possible to make a good selection of primes and to implement the
actual computation of~$\Lambda_S$ in a reasonably efficient manner, so that
a standard PC (standard as of~2008) can perform the calculations in less
than a day. We finally obtain the result we were suspecting from the beginning.

\begin{theorem}[Bugeaud, Mignotte, Siksek, Stoll, Tengely] \strut \\
   Let $x, y$ be integers satisfying
   \[ \binom{y}{2} = \binom{x}{5}\,. \]
   Then \quad
   $ x \in \{0, 1, 2, 3, 4, 5, 6, 7, 15, 19\} \,. $
\end{theorem}

A detailed description of the method (explained using the different
example equation $y^2 - y = x^5 - x$) can be found in~\cite{BMSST}.


\end{document}